\documentclass[11pt,a4paper,reqno]{amsart} 

\usepackage{amssymb,graphicx}

\usepackage{amssymb}

\newtheorem{theo}{Theorem}[section]

\newtheorem{col}[theo]{Corollary}

\newcommand{\koniec}{\begin{flushright}  $\Box $ \end{flushright}}

\newcommand{\mnote}[1]
{\protect{\stepcounter{mnotecount}}$^{\mbox{\footnotesize
$
\bullet$\themnotecount}}$ \marginpar{
\raggedright\tiny\em
$\!\!\!\!\!\!\,\bullet$\themnotecount: #1} }

\newcommand{\hook}{{\setlength{\unitlength}{11pt}   
                   \begin{picture}(.833,.8)
                   \put(.15,.08){\line(1,0){.35}}
                   \put(.5,.08){\line(0,1){.5}}
                   \end{picture}}}
\newcommand{\CP}{\mathbb{CP}}
\newcommand{\C}{\mathbb{C}}
\newcommand{\PP}{\mathbb{P}}
\newcommand{\RP}{\mathbb{RP}}
\newcommand{\R}{\mathbb{R}}

\def\p{\partial}
\def\be{\begin{equation}}

\def\ee{\end{equation}}

\def\bea{\begin{eqnarray}}
\def\eea{\end{eqnarray}}

\newcommand{\spp}{\mathbb{S}}

\begin{document}\date{January 6, 2014}
\title{Einstein--Weyl geometry, dispersionless Hirota equation and Veronese webs}
\author{Maciej Dunajski}
\address{Department of Applied Mathematics and Theoretical Physics\\ 
University of Cambridge\\ Wilberforce Road, Cambridge CB3 0WA, UK.
}\email{m.dunajski@damtp.cam.ac.uk}
\author{Wojciech Kry\'nski}
\address{Department of Applied Mathematics and Theoretical Physics\\ 
University of Cambridge\\ Wilberforce Road, Cambridge CB3 0WA, UK \\and\\
Institute of Mathematics of the  Polish Academy of Sciences\\
\'Sniadeckich 8, 00-956 Warsaw, Poland.}
\email{krynski@impan.pl}

\begin{abstract}
We exploit the correspondence between the three--dimensional Lorentzian Einstein--Weyl geometries of the hyper--CR type,
and the Veronese webs to show that the former structures are locally given in terms of solutions to the dispersionless Hirota equation.  We also demonstrate how to 
construct hyper--CR Einstein--Weyl structures by Kodaira deformations
of the flat twistor space $T\CP^1$, and how to recover the pencil of Poisson structures in five dimensions illustrating the method by an example of the Veronese web on the Heisenberg group.

\end{abstract}
\maketitle
\section{Introduction}
The notion of three--dimensional Veronese webs appearing in  the study of finite dimensional bi--Hamiltonian system is based on existence of one--parameter families of foliations of an open set of $\R^3$ 
by surfaces \cite{G,Zah,Pan}. The same structure
underlies the Einstein--Weyl geometry in 2+1 dimensions, where
the surfaces are required to be totally geodesics with respect to  some torsion--free connection compatible with a conformal structure \cite{Cartan,hitchin}. The connection between
the Veronese webs and Einstein--Weyl geometry has not so far been made, and
one purpose of this note is to show that the three--dimensional Veronese webs correspond
to a subclass of Einstein--Weyl structures which arise as symmetry reductions of
(para) hyper--Hermitian structures in signature $(2, 2)$. This class of Einstein--Weyl
structures is called hyper--CR as it admits a 
hyperboloid of (para) CR structures \cite{D1}.
In the next Section  we shall make use of this correspondence to show that
all hyper--CR Einstein--Weyl structures locally arise from solutions of the dispersionless
Hirota equation (Theorem \ref{theorem_web}). In Section \ref{section_2}
we shall elucidate the procedure of recovering a hyper--CR Einstein--Weyl structure
(and thus also a Veronese web) from the corresponding twistor space. In Section \ref{section_3} (Theorem \ref{theo_poisson}) we give an explicit 
construction of the bi--Hamiltonian structure in five dimensions from the hyper--CR Einstein--Weyl structures. 
In this construction the twistor function corresponds to the Casimir of the Poisson pencil. Finally we use the example of a Veronese web on the 
three--dimensional Heisenberg group to illustrate the procedure of recovering a solution of the dispersionless Hirota equation from a given three--parameter
family of twistor curves.
\section{Einstein--Weyl structures from dispersionless Hirota equation}
\label{section_1}
Consider the dispersionless integrable equation of Hirota type \cite{Zah,Fer,Marvan} (see also \cite{doliwa} for a discussion of its discrete version) 
\be
\label{hirota}
(b-a)\;w_x w_{yz}+a\; w_y w_{zx}-b\; w_z w_{xy}=0,
\ee
where $w=w(x, y, z)$ is the unknown function on an open set $B\subset \R^3$ and $a, b$ are non--zero constants such that
$a\neq b$. This equation arises as the Frobenius integrability condition $[L_0, L_1]=0$
for the dispersionless Lax pair of vector fields
\be
\label{lax}
L_0=\p_z-\frac{w_z}{w_x}\p_x+\lambda \;a\;\p_z, \quad L_1=\p_y-\frac{w_y}{w_x}\p_x+\lambda\; b\;\p_y.
\ee
This Lax pair is linear in the spectral parameter $\lambda$ and does not contain 
derivatives with respect to this parameter. It therefore fits into the formalism 
of \cite{D1}, and one expects equation (\ref{hirota}) to give rise to an Einstein--Weyl structure
on $B$
of hyper-CR type. To construct this structure we could find a linear
combination of $L_0, L_1$ of the form $V_1-\lambda V_2, V_2-\lambda V_3$, where
$V_i, i=1, 2, 3$ are vector fields on $B$, then read off the contravariant conformal structure
$V_2\odot V_2-V_1\odot V_3$ and try to solve a system  of differential equations
for the Einstein--Weyl one--form. We shall instead use another procedure which will
give the one--form directly (yet another, more straightforward method applicable to a broad class
of dispersionless integrable PDEs has been proposed recently \cite{FK}). Extend (\ref{lax}) to a Lax pair
of vector fields on a four--manifold $M=B\times \R$, where $\tau$ is the coordinate on $\R$
\[
{L_0}'=L_0+\p_\tau, \quad {L_1}'=L_1.
\]
This Lax pair is of the form ${L_0}'=W_1-\lambda W_2,  {L_1}'=W_3-\lambda W_4$, where
$W_1, W_2, W_3, W_4$ are linearly independent vector fields on $M$. Therefore it
defines an anti--self--dual hyper--hermitian conformal structure of neutral signature
\cite{D0,D2,zam} given in its contravariant form by\footnote{The vector fields 
$L_0', L_1'$ span an $\RP^1$ worth of 
surfaces - called $\alpha$ surfaces -- through each point of $M$. The resulting conformal structure $[g]$ is uniquely defined by demanding that these surfaces are totally isotropic and self--dual. The Frobenius integrability condition $[L_0', L_1']=0$ modulo $L_0'$ and $L_1'$ is equivalent to vanishing of the self--dual part of the Weyl tensor of $[g]$. In general the distribution $L_0', L_1'$ is defined on the $\RP^1$ bundle of real projective spinors over $M$, and the vector fields $L_0', L_1'$ contain derivatives
w.r.t the parameter $\lambda$. If no such derivatives are present, then $[g]$ contains a metric which is hyper--Hermitian \cite{D0,D2},
i. e. there exist three real anti-commuting endomorphisms $I, S, T:TM\rightarrow TM$ such that
\[
S^2=T^2=IST=-I^2={\bf 1},
\]
and any $g\in [g]$ is hermitian w.r.t the hyperboloid $a^2-b^2-c^2=1$
of complex structures  ${\bf J}=aI+bS+cT$.
} $W_1\odot W_4-W_2\odot W_3$.
Let $g$ be the corresponding covariant form.
This hyper--hermitian structure admits a non--null Killing vector $K=\p/\p \tau$, and thus
the Jones--Tod procedure \cite{JT} gives rise to an Einstein--Weyl structure on the space
of orbits $B$ of the group action generated by $K$ in $M$. This is in general given
by
\be
\label{jteq}
h=|K|^{-2}g-|K|^{-4}{\bf K}\otimes {\bf K}, \quad \omega =2|K|^{-2}*_g ({\bf K}\wedge d{\bf K}),
\ee 
where ${\bf K}=g(K, \cdot )$, $|K|^2=g(K, K)$, and $*_g$ is the Hodge operator of $g$.
This defines a Weyl structure consisting of a conformal structure
$[h]=\{ch\ | \ c\colon B\rightarrow \R^+\}$, and a torsion--free connection $D$ which is compatible with $[h]$
is the sense that 
\[
Dh=\omega\otimes h.
\]
This compatibility condition is invariant under the transformation
$h\rightarrow \phi^2 h, \;\;\omega\rightarrow \omega+2d(\ln{\phi})$, where
$\phi$ is a non--zero function on $B$.

Applying the Jones--Tod procedure (\ref{jteq})  to the Lax pair (\ref{lax}),
rescaling the resulting three--metric by $\phi^2=w_z/(w_x w_y)$, and using
(\ref{hirota}) to simplify the resulting one--form yields the Weyl structure
\begin{eqnarray}
\label{ew_final}
h&=&\frac{w_x}{w_yw_z}dx^2+\frac{a^2\;w_y}{(a-b)^2\; w_xw_z}dy^2+ \frac{b^2\;w_z}{(a-b)^2\;w_xw_y}dz^2\\
&&+2\frac{a}{(a-b)\;w_z}dxdy-2\frac{b}{(a-b)\;w_y}dxdz+2\frac{ab}{(a-b)^2\;w_x}dydz,\nonumber\\
\omega&=&-\frac{w_{xx}}{w_x}dx-\frac{w_{yy}}{w_y}dy-\frac{w_{zz}}{w_z}dz.\nonumber
\end{eqnarray}

  A Weyl structure is said to be Einstein--Weyl if the symmetrised Ricci tensor of $D$
is proportional to some metric $h\in [h]$. This conformally invariant
condition can be formulated directly as a set of non--linear PDEs on the pair
$(h, \omega)$: 
\[
R_{ij}+\frac{1}{2}\nabla_{(i}\omega_{j)}+\frac{1}{4}\omega_i\omega_j-\frac{1}{3}
\Big(R+\frac{1}{2}\nabla^k\omega_k+\frac{1}{4}\omega^k\omega_k\Big) h_{ij}=0,
\]
where $\nabla$, $R_{ij}$, and $R$ are respectively the Levi--Civita connection, the Ricci tensor and the Ricci scalar of $h$. A direct computation shows that the Weyl structure 
(\ref{ew_final}) is Einstein--Weyl iff the dispersionless Hirota equation (\ref{hirota}) holds. We must have $w_xw_yw_z\neq 0$ for the conformal structure to be non--degenerate. For example $w=x+y+z$ gives a flat Einstein--Weyl structure,
with $\omega=0$ and $h$ the Minkowski metric on $B=\R^{2, 1}$.

 
  According to Gelfand and Zakharevich \cite{G}, a Veronese web on a three--manifold $B$ is a one--parameter family of foliations of $B$ by surfaces, 
such that the normal vector fields to these surfaces through any $p\in B$ form a Veronese curve in $\PP (T_p^* B)$. The Veronese web underlying the PDE (\ref{hirota}) corresponds to a dual
 Veronese curve in $\PP (T_pB)$
\[
\lambda \longrightarrow V(\lambda)=V_1-2\lambda V_2+\lambda^2 V_3,
\]
where $\lambda\in \RP^1$, and the vector fields $V_i$ are given by
\[
V_1=
ab(a w_y\p_z-b w_z\p_y), \quad V_2=ab (w_y\p_z-w_z\p_y),
\quad
V_3=b w_y\p_z-a w_z\p_y+(a-b)\frac{w_y w_z}{w_x}\p_x.
\]
One of the seminal results of \cite{Zah} is that all Veronese webs in three dimensions are locally of this form, and thus arise from solutions to equation (\ref{hirota}).

We observe that $h(V(\lambda), V(\lambda))=0$ for all $\lambda$, where $h$ is given
by (\ref{ew_final}), so $V(\lambda)$ determines a pencil of null vectors. The vector fields
$V_1-\lambda V_2$ and $V_2-\lambda V_3$ (or equivalently the vector fields 
$L_0$ and $L_1$ given by (\ref{lax})) form and orthogonal complement of $V(\lambda)$.
For each $\lambda\in \RP^1$ they span a null surface in $B$ which is totally geodesic with respect to $D$. The existence of such totally geodesic surfaces for a triple
$(B, [h], D)$ is equivalent to the Einstein--Weyl condition \cite{Cartan}. Thus Veronese webs in three dimensions form a subset of Einstein--Weyl structures. The necessary and sufficient condition for an Einstein--Weyl structure to correspond to a Veronese web is that its  Lax pair does not contain derivatives w.r.t. $\lambda$. We have therefore established
\begin{theo}
\label{theorem_web}
There is a one to one correspondence between three--dimensional Veronese webs and Lorentzian Einstein--Weyl structures of hyper-CR type.
All hyper--CR Einstein--Weyl structures are locally of the form {\em(\ref{ew_final})}, where the function $w=w(x, y, z)$ satisfies equation
{\em(\ref{hirota})}.
\end{theo}
This Theorem, together with results of \cite{krynski} implies that there is a one to one
correspondence between hyper--CR Einstein--Weyl structures and third order ODEs such that certain two fibre preserving contact invariants vanish.

\subsubsection*{Twistor theory}
 In the real analytic category
one can make use of a complexified setting, where $w$ is a holomorphic function
on a complex domain $B_\C\subset \C^3$, and the Einstein--Weyl structure is also holomorphic.  The  twistor space of the PDE
(\ref{hirota})
is
a two--dimensional complex manifold ${\mathcal Z}$
whose points correspond to totally geodesic surfaces in $B_\C$. This gives rise to a double fibration picture
\be
\label{double_fib}
B_\C\xleftarrow{\nu} B_\C\times\CP^1\xrightarrow{\mu} {\mathcal Z},
\ee
where ${\mathcal Z}$ arises as a factor space of $B_\C\times\CP^1$ by the rank--two
distribution $\{L_0, L_1\}$ spanned by (\ref{lax}). The points in $B_\C$ 
correspond to rational curves, called twistor curves, in ${\mathcal Z}$ with self--intersection two, 
i.e. the normal bundle is $N(l_p)={\mathcal O}(2)$, where
$l_p=\mu\circ \nu^{-1}(p)$, for $p\in B_\C$. 

The special case of the hyper--CR Einstein--Weyl structures
is characterised by the existence of a holomorphic fibration of ${\mathcal Z}$ over a complex projective line \cite{D1}. 
The existence of this fibration is a consequence of the absence of vertical $\p/\p\lambda$ terms in the Lax pair (\ref{lax}).
The real manifold $B$ (and thus 
the real solutions to (\ref{hirota})) corresponds to twistor curves
invariant under an anti--holomorphic involution 
$\tau:{\mathcal Z}\rightarrow {\mathcal Z}$ which fixes an equator of each twistor curve.

\subsubsection*{Hierarchies}
Equation (\ref{hirota}) can be embedded into an 
infinite hierarchy of integrable PDEs, which
arises from the distribution
\[
L_i=\p_i-\frac{\p_i w}{\p_x w}\p_x+\lambda \;a_i\;\p_i, \quad i=0, 1, \dots
\]
and is given by $[L_i, L_j]=0$, or
\be
\label{hierarchy}
(a_i-a_j)\;\p_x w \p_i\p_jw+ a_j \p_i w \p_j\p_x w-a_i \p_j w \p _i\p_x w=0, \quad
\mbox{(no summation)}.
\ee
Here $\p_i=\p/\p x^i$ and $a_i$ are distinct constants. The function $w=w(x, x_0, x_1, \dots)$ depends on an infinite number
of independent variables, but the Cauchy data for the overdetermined system (\ref{hierarchy}) only depends on functions of two--variables which is the same functional degrees of freedom as in (\ref{hirota}).  
\section{Hyper--CR equation and Kodaira deformation theory}
\label{section_2}
In \cite{D1} it was shown that all Lorentzian hyper-CR Einstein--Weyl spaces 
on an open set $B\subset \R^3$
are locally of the form 
\be
\label{cr_metric}
h=(dY+H_X dT)^2-4(dX -H_Y dT)dT, \quad \omega=
H_{XX}dY +(H_XH_{XX}+2H_{XY})dT,
\ee
where $H=H(X, Y, T)$ satisfies the
the hyper-CR equation \cite{Dphil, pav, Fer2, ovs, manakov, BCC}
\be
\label{cr}
H_{XT}-H_{YY}+H_YH_{XX}-H_XH_{XY}=0.
\ee
Combining this result with Theorem \ref{theorem_web} suggests that
PDEs (\ref{hirota}) and (\ref{cr}) are equivalent, as there exists a local diffeomorphism of $B$ mapping (\ref{ew_final}) to (\ref{cr_metric}).
It can be shown that there is no point equivalence between these equations \cite{Fer}, but this does not rule out the existence of a  
Backlund type transformation between them which involves the second derivatives of the dependent variables (a transformation of this type
connects the first and the second Plebanski heavenly equations \cite{D2})
or a Legendre transformation of the type recently analysed by Bogdanov \cite{bog}. A procedure
recovering a solution of the Hirota equation (\ref{hirota}) from the hyper--CR equation 
(\ref{cr}) will be illustrated (at the end of Section \ref{section_3})
by the example of the Veronese web on the Heisenberg group.
 
 In this section we shall simplify the procedure of recovering the solution to (\ref{cr}) and the corresponding Einstein--Weyl structure form the twistor data \cite{D1}.  
Equation
(\ref{cr}) is equivalent to $[L_0, L_1]=0$, where 
the Lax pair is\footnote{This differs from the Lax pair of \cite{D1} by replacing $\lambda$ by $1/\lambda$.} 
\be
\label{hcr_lax}
L_0=\p_Y-\lambda(\p_T+H_Y\p_X), \quad L_1=\p_X-\lambda(\p_Y+H_X\p_X).
\ee
The ring of twistor functions spanning the kernel of this Lax pair is $(\psi, \lambda)$, 
where 
\be
\label{series_1}
\psi=T+\lambda\; Y+\lambda^2\; X+\lambda^3 \; H+\lambda ^4 \psi_4+\dots=\sum_{i=0}^\infty
\lambda^i\;\psi_i,
\ee
and the functions $\psi_i$ depend on $(X, Y, T)$ but not on $\lambda$.
The recursion relations connecting $\psi_{i+1}$ to $\psi_i$ arise from equating coefficients of $\lambda$ in
$L_0\psi=L_1\psi=0$ to zero.  The consistency conditions for these relations
imply that any of the functions $\psi_i$ satisfies
\[
(\p_X\p_T-\p_Y^2+H_Y\p_X^2-H_X\p_X\p_Y)\psi_i=0.
\]
Note that this is different than a linearisation of (\ref{cr}).
\vskip5pt

In the real analytic setup, where the twistor methods can be applied, the conformal structure $[h]$
arises from defining a null vector at $p\in B_\C$ to be a section of a normal bundle $N(l_p)$ which vanishes at one point on the twistor line 
$l_p\cong \CP^1$ in the twistor space
${\mathcal Z}$
to second order.
To define the connection $D$, define a direction at $p\in B_\C$
to be the one--dimensional space of sections on $N(l_p)$ which vanish at two
points. The one--dimensional family of curves vanishing at these two points corresponds to a geodesic in $B_\C$. In the limiting case when 
these two points coincide
this geodesic is null w.r.t $[h]$ in agreement with the definition of the Weyl structure \cite{hitchin}.
The procedure of recovering both $D$ and $[h]$ from a family of curves
in ${\mathcal Z}$ is explicit, but rather involved \cite{Mer}. This can be vastly simplified by using the local form 
(\ref{cr_metric}): given a three--parameter family of sections of
${\mathcal Z}\rightarrow \CP^1$
\[
\lambda\longrightarrow (\lambda, \psi(\lambda))
\]
choose a point on the base of the fibration ${\mathcal Z}\rightarrow \CP^1$, say $\lambda=0$, and parametrise this family by coordinates $(X, Y, T)$ where
\be
\label{parametrisation}
T=\psi|_{\lambda=0}, \quad  Y=\frac{\p \psi}{\p \lambda}|_{\lambda=0}, \quad
X=\frac{1}{2} \frac{\p^2\psi}{\p \lambda^2}|_{\lambda=0}.
\ee
Expanding the twistor function $\psi$ in a Laurent series in $\lambda$ gives
(\ref{series_1}), and we can read off $H=H(X, Y, T)$ from the coefficient of $\lambda^3$. The Einstein--Weyl structure
is now given by (\ref{cr_metric}), and we did not have to follow
the procedure of \cite{Mer}.

For example, expanding the Nil twistor function \cite{D1} 
(which depends on a constant parameter $\epsilon$, and in the limit
$\epsilon\rightarrow 0$ reduces to a section
$\lambda\rightarrow T+\lambda\; Y+\lambda^2 X$ of the
undeformed twistor space $T\CP^1$)
\be
\label{twistorf_nil}
\psi=T+\lambda\; Y -\frac{\lambda}{\epsilon}\;\ln{( 1-\lambda \epsilon X)}
\ee
in a power series in $\lambda$ and comparing with (\ref{series_1}) gives $H=\epsilon X^2/2$. 
Using (\ref{cr_metric}) we find 
the corresponding one--parameter family of deformations of the flat Einstein--Weyl structure is given by Lorentzian Nil geometry
\be
\label{Nil}
h=(dY+\epsilon XdT)^2-4dX dT, \quad \omega=\epsilon(dY+\epsilon XdT).
\ee
We shall come back to this example at the end of Section \ref{section_3}.
\subsubsection*{Deformation theory}
Let $[\psi: \pi_0: \pi_1]$ be homogeneous holomorphic coordinates
on the twistor space ${\mathcal Z}=T\CP^1$ corresponding to the flat
Einstein--Weyl structure with $H=0$. Here $[\pi_0: \pi_1]$ are homogeneous coordinates
on the base $\CP^1$ and $\lambda=\pi_0/\pi_1$.
Cover ${\mathcal Z}$
by two open sets ${\mathcal{U}}$ and $\widetilde{{\mathcal{U}}}$ with $\pi_1\neq 0$ on ${\mathcal{U}}$ and 
$\pi_0\neq 0$ on $\widetilde{{\mathcal{U}}}$. 
Following the steps of Penrose's non--linear graviton construction
\cite{Pen} we construct twistor spaces corresponding to non--trivial
Einstein--Weyl spaces by deforming the patching relation between ${\mathcal{U}}$
and $\widetilde{{\mathcal{U}}}$. The Kodaira theorem \cite{Kod} guarantees that
the deformations preserve the three--parameter family  of curves and the deformed twistor space still gives rise to a three--dimensional manifold
$(B_\C, [h], D)$ with Einstein--Weyl structure. 

The hyper--CR deformations preserve the fibration of ${\mathcal Z}$ over $\CP^1$. Thus we deform the patching relations
such that $\pi_0/\pi_1$ is preserved on the overlap. The infinitesimal deformations
of this kind are generated by elements of $H^1({\mathcal Z}, 
{\bf \Theta})$, where ${\bf \Theta}$ denotes a sheaf of germs of holomorphic vector fields. Let a vector
field
\[
{\mathcal Y}=f\frac{\p}{\p \psi}+g\Big(\pi_0\frac{\p}{\p \pi_0}+\pi_1\frac{\p}{\p \pi_1}\Big)
\]
defined on an overlap of ${\mathcal{U}}$ and $\widetilde{{\mathcal{U}}}$ define a class
in $H^1({\mathcal Z}, 
{\bf \Theta})$. Here the functions $f$ and $g$ are holomorphic and homogeneous of degree two and zero respectively on the overlap.
The finite deformation $\tilde{\psi}=\tilde{\psi}(\psi, \lambda, \epsilon)$ is given by integrating
\[
\frac{d \psi}{d\epsilon}=f, \quad \frac{d\pi_0}{d\epsilon}=g\pi_0,
\quad \frac{d\pi_1}{d\epsilon}=g\pi_1.
\] 
It is possible that the assumption about real--analyticity can be dropped by employing the techniques of
\cite{Nak} or some modification of the dispersionless $\bar{\p}$ approach
\cite{BK}, but that remains to be seen.
\section{Bi-Hamiltonian structure}
\label{section_3}
Recall that a Poisson structure on a manifold $U$ is a bivector $P\in\Lambda^2(TU)$ such that the associated
Poisson bracket $\{f, g\}:=P\hook (df\wedge dg)$ satisfies the Jacobi identity
\[
\{f,\{ g, h\}\}+
\{h,\{ f, g\}\}+\{g,\{ h, f\}\}=0
\]
for all smooth functions $f, g, h$ on $U$.
If we regard $U\subset \R^{n}$ as an open set, and adapt local coordinates such that
\[
P=\sum_{\alpha, \beta=1}^n P^{\alpha\beta}\frac{\p}{\p x^{\alpha}}\wedge \frac{\p}{\p x^{\beta}},
\]
then the  Poisson bracket takes the form
\[
\{f, g\}= \sum_{\alpha, \beta=1}^n P^{\alpha\beta}\frac{\p f}{\p x^{\alpha}} \frac{\p  g}{\p x^{\beta}}.
\]
If $n$ is odd, the Poisson structure necessarily admits at least one Casimir, i.e. a non--constant function $C$ such that $\{C, f\}=0$ for all functions $f$ on $U$.

A bi-Hamiltonian structure   $P^{\lambda}$ on $U$ consists of two Poisson structures $P_0$ and $P_1$ such that 
the pencil of Poisson brackets 
\[
\{f, g\}_\lambda:=\{f, g\}_0+\lambda \{f, g\}_1
\]
satisfies the Jacobi identity for all values of the parameter $\lambda$  (see e.g. \cite{magri,olver,blaszak}).

 We shall assume that $n=5$ and 
$P^\lambda$ has exactly one Casimir function $C^\lambda$ for each $\lambda$.
Thus, regarding the $P^{\lambda}\in \Gamma(\Lambda^2(TU))$ as a map from
$T^*U$ to $TU$, we find that this map has a one--dimensional kernel.
Therefore, for each $\lambda\in \RP^1$, the manifold $U$ is foliated by four--dimensional symplectic leaves of $P^{\lambda}$, i. e. 
there exists a rank four integrable distribution ${\mathcal D}_\lambda \subset TU$ such that 
${\mathcal D}_\lambda= {\mathrm{ker}}(d C^\lambda)$. This distribution is annihilated by a $\lambda$--dependent one--form
\[
{\bf e}(\lambda):=(P^{\lambda}\wedge P^{\lambda})\hook\Omega,
\]
where $\Omega$ is some fixed volume form on $U$, and $\hook$ denotes its
contraction with a section of $\Lambda^{4}(TU)$. 
Therefore  ${\bf e}(\lambda)$ is quadratic in $\lambda$, and we can write
\be
\label{lambda_form}
{\bf e}(\lambda)=e^3+\lambda e^2+\lambda^2 e^1
\ee
for some one--forms $e^i$. The Pfaff theorem implies that
$dC^{\lambda}\in {\mathrm{span}}({\bf e}(\lambda)).$

The intersection of all distributions ${\mathcal D}_\lambda$
as $\lambda\in \RP^1$ varies is an integrable rank--two distribution which we shall call 
${\mathcal D}$. Thus 
\[{\mathcal D}=\ker e^1\cap\ker e^2\cap\ker e^3=\bigcap_\lambda \ker
({\bf e}(\lambda))\subset TU\]
which is integrable, as an intersection of tangent bundles to the 
symplectic leaves of the Casimir $C^{\lambda}$.
Therefore
${\mathcal D}\subset {\mathcal D}_\lambda\subset TU$. We shall assume that all these distributions have constant ranks.
The two dimensional distribution ${\mathcal D}_{\mathcal Z}:={\mathcal D}_\lambda/{\mathcal D}$
is defined on the three dimensional quotient $B:=U/{\mathcal D}$. This distribution is integrable,
and it defines the structure of a Veronese web on $B$. This web is a projection of symplectic leaves of $P^{\lambda}$ to $B$. 

This is essentially the construction of 
\cite{G}. In \cite{G} an inverse construction of a bi-Hamiltonian structure from a Veronese web is
also given. This construction does not appear to be explicit, but it has been shown \cite{tur} that a bi-Hamiltonian structure on $U$ can be recovered form a Veronese web on $B$. 

We shall now
give a simple algorithm for recovering the bi-Hamiltonian structure from any hyperCR Einstein--Weyl
structure which, by Theorem \ref{theorem_web}, is equivalent to a Veronese web. In our procedure
the distribution ${\mathcal D}_{\mathcal Z}$ is identified with the twistor distribution (\ref{lax})
given by the span of the Lax pair $L_0$ and $L_1$ which (in the complexified setting) underlies the double fibration picture (\ref{double_fib}). The one 
forms $e^i$ in (\ref{lambda_form}) define the conformal structure
\[
h=e^2\otimes e^2-2(e^1\otimes e^3+e^3\otimes e^1)
\]
and the Frobenius condition ${\bf e}\wedge d{\bf e}=0$ is equivalent to the hyper--CR Einstein--Weyl equations
\cite{D1}.

Before formulating the next Theorem
recall  that the conformal metric of signature $(2, 1)$
on a three--dimensional manifold $B$ is equivalent to the existence of an isomorphism
\be
\label{para_con}
TB \cong \spp\odot \spp,
\ee
where $\spp$ is a rank two real symplectic vector bundle over $B$. In concrete terms (\ref{para_con}) identifies vectors with two by two symmetric matrices (or equivalently, with
homogeneous quadratic polynomials in two variables). The conformal structure is then defined by
declaring a vector to be null iff the corresponding matrix has rank one (or equivalently, if the 
corresponding homogeneous polynomial has  a repeated root). The isomorphism (\ref{para_con})
is a particular example of a $GL(2, \R)$ structure: the group  $GL(2, \R)$ acts
on the fibres of $\spp$, and induces Lorentz rotations and conformal rescalings on vectors in $B$
as  $\mathfrak{gl}(2, \R)\cong \mathfrak{so}({2, 1})\oplus \R$.

Let $(p_0, p_1)$ be local coordinates on the fibres of the dual vector bundle $\spp^*$.
\begin{theo}
\label{theo_poisson}
Let $(B, [h], D)$ be a hyper--CR Einstein--Weyl structure and let $(L_0, L_1)$ be the 
commuting twistor 
distribution (Lax pair) spanning a pencil of null and totally geodesic surfaces
through each point $p\in B$. Then
\be
\label{poisson}
P^\lambda =L_0\wedge\frac{\p}{\p p_0}+L_1\wedge\frac{\p}{\p p_1}.
\ee
is a bi--Hamiltonian structure on the five--dimensional manifold $\spp^*$, where 
$\spp^*$ is a real rank two bundle over $B$ such that (\ref{para_con}) is the canonical isomorphism given
by the conformal structure $[h]$ on $B$.
\end{theo}
{\bf Proof.} Consider three linearly independent vector fields
$(V_1, V_2, V_3)$ on $B$ such that the conformal structure in the Einstein--Weyl structure is represented by the contravariant form
$V_2\odot V_2-V_1\odot V_3$.
Let $\textit{}\zeta_\lambda\subset B$ be a two--dimensional totally geodesic surface through
$p\in B$ corresponding to $\lambda\in\RP^1$
and let $(p_0, p_1)$ be coordinates on the fibres of
the cotangent bundle $T^*\zeta_\lambda$. The tangent space
$T_p\zeta_\lambda$ is spanned by the twistor distribution
$L_0=V_1-\lambda V_2, L_1=V_2-\lambda V_3$ (any twistor distribution underlying a hyper-CR   
Einstein--Weyl space can be put in this form, by taking a linear combination of the spanning vector with $\lambda$--independent coefficients). 
The vector fields $L_0$ and $L_1$ commute so the canonical 
Poisson structure on $T^*\zeta_\lambda$ is
\[
L_0\wedge\frac{\p}{\p p_0}+L_1\wedge\frac{\p}{\p p_1}.
\]
We now extend this to a bi--Hamiltonian structure on 
$B\times \R^2$. To do it we need to canonically identify the tangent
spaces to different null totally geodesic  surfaces in $B$ and show that the union of cotangent bundles $\bigcup_{\lambda\in \RP^1}T^*\zeta_\lambda$ is isomorphic to $\spp^*$. To establish the latter,
observe that  under the isomorphism (\ref{para_con}) vectors 
tangent to $\zeta_\lambda$ correspond
to quadratic polynomials with one of the roots equal to $\lambda$: 
represent
the basis $(V_1, V_2, V_3)$ of $T_pB$ as a symmetric matrix of vector fields
$V_{AB}$ with $A ,B=0, 1$ and $V_{00}=V_1, V_{01}=V_{10}=V_2, V_{11}=V_3$, and the twistor distribution is $L_A=\pi^BV_{AB}$ where $\pi^B=[\pi^0, \pi^1]$ are homogeneous coordinates on 
$\RP^1$ such that $\pi^1/\pi^0=-\lambda$. Thus any vector tangent to $\zeta_\lambda$ is
$W=\mu^AL_A$ for some $\mu^A=[\mu^0, \mu^1]$. The two by two matrices corresponding  to such vectors
are $W=\mu\odot\pi$,
with $\pi^A$ fixed and $\mu^A$ varying. Equivalently vectors tangent to 
$\zeta_\lambda$ correspond to homogeneous polynomials with common root $[\pi^0, \pi^1]$. Dividing out by the factor containing this root yields
linear homogeneous polynomials, i.e. an isomorphism $T\zeta_\lambda\cong 
\spp|_{\zeta_\lambda}$.
To obtain the canonical identification of all tangent spaces $T_p\zeta_\lambda$ as $\lambda$ varies we simply evaluate the twistor distribution at different values of $\lambda$: Let
$\zeta_1$ and $\zeta_2$ be totally geodesic null surfaces corresponding to $\lambda_1$ and
$\lambda_2$ respectively. 
The isomorphism $\phi_{12}:T\zeta_1\rightarrow T\zeta_2$ is defined on the basis vectors by
\[
\phi_{12}(L_0|_{\lambda_1})=L_0|_{\lambda_2}, \quad \phi_{12}(L_1|_{\lambda_1})=L_1|_{\lambda_2}.
\]
Therefore the fibre coordinates $(p_0, p_1)$ can be used unambiguously for all $\lambda$
and the Poisson structure on $T^*\zeta_\lambda$ extends to a Poisson pencil\footnote{This agrees with the Gelfand--Zakharevich construction \cite{G}, where
the existence  of the isomorphism $T^*B=\spp^*\odot\spp^*$ is postulated. 
The authors 
consider
\[
T^*B\otimes \spp^*=(\spp^*\odot\spp^*\odot\spp^*) \oplus (\spp^*\otimes\Lambda^2(\spp^*)).
\]
The 5-manifold supporting the Poisson pencil arises as the second factor in this decomposition, i.e. the total space of $\spp^*\rightarrow B$ as the line bundle $\Lambda^2(\spp^*)$ is trivialised.} (\ref{poisson})
on $\bigcup_{\lambda\in\RP^1} {\spp^*}|_{\zeta_\lambda}=\spp^*$.

For any Poisson structure given by a bivector $P$ the Jacobi identity reduces to
\[
\sum_{\delta=1}^n P^{\delta\gamma}\frac{\p P^{\alpha\beta}}{\p x^{\delta}} + 
P^{\delta\beta}\frac{\p P^{\gamma\alpha}}{\p x^{\delta}}+
P^{\delta\alpha}\frac{\p P^{\beta\gamma}}{\p x^{\delta}}=0 \quad\mbox{for all}\quad \alpha, \beta, \gamma=1, \dots n.
\]
This set of conditions holds, with $n=5$, 
for $P^{\lambda}=P=P_0-\lambda P_1$ 
regardless of a particular form of the Lax pair $(L_0, L_1)$, as long as
$[L_0, L_1]=0$ because in this case the Frobenius theorem can be used to put
the Poisson bi-vector in the canonical form.

If the vector fields of the Lax pair are linear in the parameter $\lambda$ - which is the case for Lax pairs underlying hyperCR Einstein--Weyl structures - then the bivector (\ref{poisson}) defines a 
pencil of compatible Poisson structures.
\koniec
In the construction of the Poisson pencil (\ref{poisson})  the Casimir $C^{\lambda}$ is a twistor function which we called $\psi$ in formula (\ref{series_1}), and
\[
{\mathcal D}_{\mathcal Z}=\mbox{span}(L_0, L_1), \quad {\mathcal D}=\mbox{span}(\p/\p p_0, \p/\p p_1).
\]
Using the explicit form of the twistor distribution (\ref{lax}) yields the following
\begin{col}
Let $w=w(x, y, z)$ be a function on $B$ with continuous second derivatives.
The Poisson brackets 
\[
P_0= \left(\begin{array}{ccccc}
0&0&0&-w_z/w_x&-w_y/w_x\\
0&0&0&0&1\\
0&0&0&1&0\\
w_z/w_x&0&-1&0&0\\
w_y/w_x&-1&0&0&0
\end{array}
\right) , \quad P_1=
\left(\begin{array}{ccccc}
0&0&0&0&0\\
0&0&0&0&b\\
0&0&0&a&0\\
0&0&-a&0&0\\
0&-b&0&0&0
\end{array}
\right).
\]
define a compatible Poisson pencil 
 if and only if the function $w=w(x, y, z)$ satisfies the dispersionless Hirota equation
{\em(}\ref{hirota}{\em)}.
\end{col}
To obtain a bi-Hamiltonian system,
the Hamiltonian function needs to be selected as one of the coefficients
in the expansion of the Casimir $C^{\lambda}$ in the powers of $\lambda$. 
It follows from the general theory that all  coefficients of this expansion are in involution with respect to both Poisson structures $P_0$ and $P_1$. 
Thus the coordinate functions on $B$ are first integrals of the Hamiltonian flow, and the integral curves of this flow are vertical to $B$. They are of the form $p_k(t)=\omega_k(x,y,z)t+\rho_k(x, y, z)$, where $k=0,1$, and $\omega_k, \rho_k$ are known functions on $B$ whose form depends on the choice of the Hamiltonian. 
This simple form of solution curves reflects the fact that the Gelfand--Zakharevich construction of Veronese webs puts the system in the action--angle coordinates.
\subsubsection*{{\bf Example. The Heisenberg group}}
We have already shown that the Heisenberg group carries an Einstein--Weyl structure
of Lorentzian signature given by (\ref{Nil}). The twistor distribution
(\ref{hcr_lax}) defining the Veronese web is
\be
\label{twistor_nil}
L_0=\p_Y-\lambda\p_T,\quad L_1=\p_X-\lambda(\p_Y+\epsilon X\p_X).
\ee
The metric $h$ defining the conformal  structure is left--invariant with the Killing vectors generated by the right invariant vector fields
\[
R_Y=\p_Y,\quad R_T=\p_T, \quad R_X=\p_X-\epsilon T\p_Y
\]
forming the Heisenberg Lie algebra
\[
[R_X, R_T]=\epsilon\;R_Y, \quad [R_X, R_Y]=0, \quad [R_T, R_Y]=0.
\]
The Einstein--Weyl structure is only invariant under the abelian subalgebra spanned by 
$R_Y$ and $R_T$ as the Lie derivative of the twistor distribution along $R_X$ is non--zero.
The Poisson pencil is also not invariant under the Heisenberg action, which is best seen by writing it in terms of the left--invariant vector fields (which commute with the right-invariant vector fields) as
\[
P^\lambda=(L_Y\wedge\p_{p_0}+L_X\wedge\partial_{p_1})
-\lambda(
(L_T+\epsilon XL_Y)\wedge\partial_{p_0}+(L_Y+\epsilon XL_X)\wedge\partial_{p_1})
\]
and ${\mathcal{L}}_{R_X} P^{\lambda}=\epsilon\lambda \; P^0\neq 0$. We note that 
there is one Poisson bracket in the pencil (corresponding to $\lambda=0$) which is left--invariant.

We finish the discussion of this example by presenting the solution to the dispersionless
Hirota equation (\ref{hirota}) corresponding to the Heisenberg Veronese web (\ref{twistor_nil}). The procedure we 
shall use elucidates two different parametrisations of twistor curves leading to 
equations (\ref{cr})
and (\ref{hirota}) respectively. In case of the hyper--CR equation we choose a point, say $\lambda=0$, on the base of the fibration 
${\mathcal Z}\rightarrow \CP^1$ (here we again work in the complexified 
setting), and
parametrise the curve $l_p$ given by $\lambda\rightarrow (\lambda, \psi(\lambda))$
corresponding to $p\in B$ by the position of its intersection with the fibre over $\lambda=0$ and the first two jet coordinates as in (\ref{parametrisation}). The third
jet then defines the dependent variable $H$
in equation (\ref{cr}).  This can be seen from formula (\ref{twistorf_nil}).

In case of the dispersionless Hirota equation we instead choose four distinct points
$\lambda_1, \lambda_2, \lambda_3$ and $\lambda_4$ on the base. The intersections of the curve
with the three fibres over $\lambda_1, \lambda_2, \lambda_3$ give, up to a reparametrisation,
the coordinates $(x, y, z)$. The intersection with the fibre over $\lambda_4$ then defines
the dependent variable $w$. We can use the M\"obius transformation to map $(\lambda_1, \lambda_2, \lambda_3)$ to $(0, \infty, 1)$. Thus we expect the coordinate $T$ in the hyper--CR equation to map directly to
coordinate $x$ in the dispersionless Hirota equation. 
This is the twistorial interpretation of the original procedure of Zakharevich
where the two--dimensional distribution ${\mathcal D}_{\mathcal Z}$
defining the  Veronese web is
\begin{eqnarray*}
{\mathcal D}_{\mathcal Z}(\lambda)&=&\ker\Big((\lambda_1-\lambda_4)(\lambda-\lambda_2)(\lambda-\lambda_3)w_xdx\\
&+&
(\lambda_2-\lambda_4)(\lambda-\lambda_1)(\lambda-\lambda_3)w_ydy\\
&+&
(\lambda_3-\lambda_4)(\lambda-\lambda_1)(\lambda-\lambda_2)w_zdz\Big).
\end{eqnarray*}
Comparing this with ${\mathcal D}_{\mathcal Z}(\lambda)=\mbox{span}(L_0, L_1)$ given by the twistor distribution (\ref{twistor_nil})
and using the M\"obius freedom described above yields
${\mathcal D}_{\mathcal Z}(0)=\mbox{ker\;}(dx), {\mathcal D}_{\mathcal Z}(1)=\mbox{ker\;}(dz), {\mathcal D}_{\mathcal Z}(\infty)=\mbox{ker\;}(dy)$, or equivalently
\[
d\psi|_{\lambda=0}\in\mbox{span}(dx), \quad
d\psi|_{\lambda=1}\in\mbox{span}(dz), \quad
d\psi|_{\lambda=\infty}\in\mbox{span}(dy), \quad
\]
where $\psi$ is the twistor function (\ref{twistorf_nil}).
These formulae define the coordinates up to reparametrisations
$(x, y, z)\rightarrow (\hat{x}(x), \hat{y}(y), \hat{z}(z))$, and we choose
\[
x=T, \quad y=Xe^{-\epsilon\;Y}, \quad z=(1-\epsilon\;X)e^{-\epsilon(T+Y)}.
\]
The solution $w(x, y, z)$ of the dispersionless Hirota equation (\ref{hirota}) can now be read--off from ${\mathcal D}_{\mathcal Z}(\lambda_4)=\mbox{ker\;}(dw)$ which yields 
$w=(1-\epsilon\lambda_4 X)e^{-\epsilon({\lambda_4}^{-1}T+Y)}$, or after transforming to the new variables and reabsorbing some constants by rescaling $(x, y, z)$,
\[
w=y e^{ax}+z e^{bx},
\]
where $\lambda_4=1-b/a$.
\subsubsection*{Acknowledgements} The work of Wojciech Kry\'nski has been partially supported by the Polish National Science Centre grant ST1/03902.

\end{document}